\documentclass{conm-p-l}
\usepackage{amsfonts}

\def\Cal{\mathcal}

\def\B{{\Cal B}}
\def\H{{\Cal H}}

\def\M{{\Cal M}}
\def\P{{\Cal P}}

\def\S{{\Cal S}}
\def\F{{\Cal F}}

\def\I{{\Cal I}}

\def\W{{\Cal W}}

\def\d{\partial}
\def\tr{{\hbox{\rm tr}}}

\def\Ma{\frM_{n,m}}

\def\W{\mathcal{W}}

\def\f0{f_0}
\def\Fc0{\varphi_0}
\def\rn{\bbr^n}

\def\I_k {I_{-}^{k/2}}
\def\I+k {I_{+}^{k/2}}

\def\bbr{{\Bbb R}}

\def\tr{{\hbox{\rm tr}}}

\def\cos{{\hbox{\rm cos}}}
\def\det{{\hbox{\rm det}}}

\def\part{\partial}
\def\intl{\int\limits}
\def\b{\beta}

\def\Gam{\Gamma}

\def\a{\alpha}

\def\Del{\Delta}
\def\del{\delta}
\def\vp{\varphi}

\def\gam{\gamma}
\def\Lam{\Lambda}
\def\sig{\sigma}
\def\lam{\lambda}

\def\e{\varepsilon}
\def\t{\tau}

\font\frak=eufm10

\def\fr#1{\hbox{\frak #1}}

\def\frM{\fr{M}}

\def\gma{\Gamma_m(\a)}

\def\cos{{\hbox{\rm cos}}}
\def\det{{\hbox{\rm det}}}

\def\p{\P_m}
\def\gm{\Gamma_m}
\def\tr{{\hbox{\rm tr}}}
\def\part{\partial}
\def\intl{\int\limits}
\def\b{\beta}

\def\Gam{\Gamma}

\def\a{\alpha}
\def\cpm{\overline\P_m}

 \newtheorem{theorem}{Theorem}[section]
 \newtheorem{lemma}[theorem]{Lemma}
 \theoremstyle{definition}
 \newtheorem{definition}[theorem]{Definition}
\newtheorem{example}[theorem]{Example}

\theoremstyle{remark}
\newtheorem{remark}[theorem]{Remark}

 \numberwithin{equation}{section}

\newcommand{\be}{\begin{equation}}
\newcommand{\ee}{\end{equation}}

\newcommand{\bea}{\begin{eqnarray}}
\newcommand{\eea}{\end{eqnarray}}
\newcommand{\Bea}{\begin{eqnarray*}}
\newcommand{\Eea}{\end{eqnarray*}}

\begin{document}

\title[Composite Wavelet Transforms]{ Composite Wavelet Transforms: Applications and Perspectives}

\author{Ilham A. Aliev}
\address{Department of Mathematics, Akdeniz University, 07058 Antalya TURKEY}
\email{ialiev@akdeniz.edu.tr}
\author{Boris Rubin}
\address{Department of Mathematics, Louisiana State University, Baton Rouge,
Louisiana 70803}
\email{borisr@math.lsu.edu}
\author{Sinem Sezer}
\address{Faculty of Education, Akdeniz University, 07058 Antalya TURKEY}
\email{sinemsezer@akdeniz.edu.tr}
\author{Simten B. Uyhan}
\address{Department of Mathematics, Akdeniz University, 07058 Antalya TURKEY}
\email{simten@akdeniz.edu.tr}
\thanks{The research was supported by the Scientific Research Project
Administration Unit of the Akdeniz University (Turkey) and TUBITAK
(Turkey). The second author  was also supported  by the  NSF grants
EPS-0346411 (Louisiana Board of Regents) and DMS-0556157.}

\renewcommand{\subjclassname}{%
    \textup{2000} Mathematics Subject Classification}
\subjclass[2000]{42C40, 44A12, 47G10.}


\keywords{Wavelet transforms, potentials, semigroups, generalized
translation, Radon transforms,  inversion formulas, matrix spaces.}

\begin{abstract}
We introduce a new concept of the so-called {\it composite wavelet
transforms}. These  transforms are  generated by two components,
namely, a kernel function and a wavelet function (or a measure). The
composite wavelet transforms and the relevant Calder\'{o}n-type
reproducing formulas constitute a unified approach to explicit
inversion of the Riesz, Bessel, Flett, parabolic and some other
 operators of the potential type generated by ordinary (Euclidean) and generalized
(Bessel) translations. This approach is exhibited in the paper.
 Another concern  is application of the  composite
wavelet transforms to explicit inversion of the k-plane Radon
transform on $\bbr^n$. We also discuss in detail  a series of open
problems arising in wavelet analysis of $L_p$-functions of matrix
argument.

\end{abstract}

\maketitle

\centerline{Contents} \centerline{} 1. Introduction.

2. Composite wavelet transforms for  dilated kernels.

3. Wavelet transforms associated to one-parametric semigroups and
inversion ${}\qquad {}\quad$ of potentials.

4. Wavelet transforms with the generalized  translation operator.

5. Beta-semigroups.

6. Parabolic wavelet transforms.

7. Some applications to inversion of the $k$-plane Radon transform.

8. Higher-rank composite wavelet transforms and open problems.

References.

\section{Introduction}

Continuous wavelet transforms
\begin{equation*}
\mathcal{W}f(x,t)=t^{-n}\int_{\bbr^n}f(y)\, w \left (\frac{x-y}{t}\right )\, dy, \qquad
x\in \mathbb{R}%
^{n},\ \ \ t>0,
\end{equation*}%
where $w$ is an integrable radial function satisfying
$\int_{\rn}w(x)dx=0$, have proved to be a powerful tool in analysis
and applications. There is a vast literature on this subject (see,
e.g., \cite {Da}, \cite{16}, \cite{20}, just for few). Owing to the
formula \bigskip
\begin{equation}
\int_{0}^{\infty }\mathcal{W}f(x,t)\text{ }\frac{dt}{t^{1+\alpha }}%
=c_{\alpha ,w}(-\Delta )^{\alpha /2}f(x),\qquad \alpha \in \mathbb{C},\quad
\Delta =\sum\limits_{k=1}^{n}\frac{\partial ^{2}}{\partial
x_{k}^{2}},  \label{1.2}
\end{equation}
that can be given  precise meaning, continuous wavelet transforms
enable us to resolve a variety of problems dealing with powers of
differential operators. Such problems arise, e.g., in potential
theory, fractional calculus, and integral geometry; see, \cite{16},
\cite{21}-\cite{27}, \cite{32}. Dealing with functions of several
variables, it is always tempting to reduce the dimension of the
domain of the wavelet function $w$ and  find new tools to gain extra
flexibility. This is actually a motivation for our article.

We introduce a new concept of the so-called {\it composite wavelet
transforms}. Loosely speaking, this is a class of wavelet-like
transforms generated by two components, namely, a kernel function
and a wavelet. Both are in our disposal. The first one depends on as
many variables as we need for  our problem. The second component,
which is a wavelet function (or a measure), depends only on one
variable. Such transforms are usually associated with one-parametric
semigroups, like Poisson, Gauss-Weierstrass, or metaharmonic  ones,
and can be implemented to obtain explicit inversion formulas for
diverse operators of the potential type and fractional integrals.
These arise in integral geometry in a canonical way; see, e.g.,
\cite{15, 22,
 26, R9}.

In the present article we study different types of composite wavelet
transforms in the framework of the $L_p$-theory and the relevant
Fourier and Fourier-Bessel harmonic analysis. The main focus is
 reproducing formulas of Calder\'{o}n's type and explicit inversion of Riesz, Bessel,
 Flett, parabolic, and some other  potentials. Apart of a brief review
 of recent developments in the area, the paper contains a
 series of
 new results. These include  wavelet transforms for
 dilated kernels and wavelet transforms generated by Beta-semigroups associated
 to multiplication by  $\exp(-t|\xi|^{\b}), \; \b>0$,  in terms of the Fourier transform.
 Such semigroups arise in the context of stable random processes in probability
 and enjoy a number of remarkable properties \cite{Ko}, \cite{La}.
 Special emphasis is made on
 detailed discussion of  open
 problems arising in wavelet analysis of functions of matrix
 argument. Important results for $L_2$-functions in this ``higher-rank"  set-up were obtained in \cite{OOR} using the Fourier transform
 technique. The $L_p$-case for $p\neq 2$ is still mysterious. The main
 difficulties are related to correct definition and handling  of admissible wavelet functions on the cone of positive definite symmetric matrices.

 The paper is organized according to the Contents
 presented above.

\section{Composite Wavelet Transforms for  Dilated Kernels}

\subsection{Preliminaries}
 Let $L_{p}\equiv L_{p}(\mathbb{R}^{n}), \; 1\le p<\infty,$ be the standard space of functions with
the norm
\begin{equation*}
\left\| f\right\| _{p}=\Big( \int_{\mathbb{R}^{n}}\left| f(x)\right|
^{p}dx\Big )^{1/p}<\infty.
\end{equation*}
For technical reasons, the notation $L_{\infty }$ will be used for
the space $C_{0}\equiv C_{0}(\mathbb{R}^{n})$  of all continuous
functions on $ \mathbb{R}^{n}$ vanishing at infinity. The Fourier
transform of a function $f$ on $\mathbb{R}^{n}$ is defined by
\begin{equation*}
Ff(\xi)=\int_{\mathbb{R}^{n}}f(x)\, e^{ix\cdot \xi }\,dx, \qquad x
\cdot \xi=x_{1}\xi _{1}+\cdots +x_{n}\xi _{n}.
\end{equation*}
For $ 0\le a<b\le\infty$, we  write $\int_a^b f(\eta)d\mu (\eta)$ to
denote the integral of the form $\int_{[a,b)} f(\eta)d\mu (\eta)$.

\noindent\begin{definition}\label{d1}  Let $q$ be a measurable
function  on $\bbr^n$ satisfying the following conditions:

(a) $q\in L_1\cap L_r$ for some $r>1$;

(b) the least radial decreasing majorant  of $q$ is integrable, i.e.
$$\tilde q (x)=\sup_{|y|>|x|} |q(y)| \in L_1;$$

(c)$\qquad \int_{\bbr^n} q(x)\, dx =1.$

\noindent We denote \be\label {qu} q_t (x)=t^{-n} q(x/t), \qquad Q_t
f(x)=(f*q_t)(x), \qquad t>0,\ee and set \be \label{cwtr}Wf(x,t)=
\int_0^\infty Q_{t\eta} f(x)\,d\mu (\eta),\ee where $\mu$ is a
finite Borel measure on $[0,\infty)$. If $\mu$ is {\it a wavelet
measure} (i.e., $\mu$
 has a certain number of vanishing moments and obeys suitable decay conditions)
 then (\ref{cwtr}) will be called the {\it composite wavelet transform} of
 $f$. The function $q$ will be called
 a {\it kernel function}  and $Q_t$ a {\it kernel operator}
 of the composite transform $W$.
\end {definition}

The integral (\ref{cwtr}) is well-defined for any function $f \in
L_p$, and
$$
||Wf(\cdot,t)||_p\le ||\mu || \, ||q||_{1} \, ||f||_p,
$$
where $||\mu ||=\int_{[0,\infty)}d|\mu| (\eta)$. We will also
consider a more general weighted transform
 \be
\label{cwtrw}W_af(x,t)= \int_0^\infty  Q_{t\eta}
f(x)\,e^{-at\eta}\,d\mu (\eta),\ee where $a\geq 0$ is a fixed
parameter.

The kernel function $q$, the wavelet measure $\mu$, and the
parameter $a \geq 0$ are in our disposal. This feature makes the new
transform  convenient in applications.

\subsection{Calder\`on's identity} An analog of Calder\'on's reproducing formula  for $W_af$ is
 given by the following theorem.
\begin{theorem}\label{teo:34}
 Let $\mu$ be  a finite Borel measure on $[0,\infty)$ satisfying
\be\label{eq:32}
 \mu([0,\infty))=0 \quad  \text{and} \quad  \int_0^\infty |\log \eta|\, d|\mu|(\eta) <\infty.\ee
If $f\in L_p, \; 1 \le p \le \infty$\footnote{We remind that
$L_\infty$ is interpreted as the space $C_0$ with the uniform
convergence.}, and
$$c_\mu=\int_0^\infty \log \frac{1}{\eta} \,\, d\mu (\eta),$$ then \be\label{eq33}\int_0^\infty
W_af(x,t)\frac{dt}{t}\equiv \lim\limits_{\e \rightarrow
0}\int_\e^\infty W_af(x,t)\frac{dt}{t}=c_\mu f(x) \ee where the
limit exists in the $L_p$-norm and pointwise for almost all $x$.  If
$f\in C_0 $, this limit is uniform on $\rn$.
\end{theorem}
\begin{proof}
 Consider
the truncated integral
\begin{equation}\label{eq:38}
I_\e f(x)=\int_\e^\infty  W_af(x,t)\frac{dt}{t}, \qquad \e >0.
\end{equation}
Our aim is to represent it in the form
\begin{equation}\label{eq:3100}
I_\e f(x)= \int_{0}^{\infty}Q_{\e s} f(x) \,e^{-a\e s}\, k(s) \, ds
\end{equation}
where \be\label{ka} k\in L_1 (0, \infty) \qquad \text{\rm and}
\qquad \int_0^\infty k(s) ds=c_\mu. \ee Once (\ref{eq:3100}) is
established, all the rest follows from properties (a)-(c) in
Definition \ref{d1} according to the standard machinery of
approximation to the identity; see \cite {St}.

Equality (\ref{eq:3100}) can be formally obtained by changing the
order of integration, namely,
\begin{eqnarray}
I_\e f(x) &=&\int_0^\infty  d\mu(\eta) \int_{\e }^\infty Q_{t\eta}
f(x)\,e^{-at\eta}\,\frac{dt}{t}\nonumber
\\&=& \int_{0}^{\infty}d\mu(\eta) \int_{\eta}^\infty Q_{\e s}
f(x)\,e^{-a\e s}\,\frac{ds}{s}\nonumber
\\&=&\int_{0}^{\infty}Q_{\e s}
f(x)\,e^{-a\e s} k(s)\, ds, \qquad
k(s)=s^{-1}\int_{0}^{s}d\mu(\eta). \nonumber
\end{eqnarray}
Furthermore, since $\mu([0,\infty))=0$, then \bea \int_0^\infty
|k(s)| ds&=&\int_0^1 \Big |\int_0^s d\mu(\eta)\Big
|\frac{ds}{s}+\int_1^{\infty} \Big |\int_s^{\infty} d\mu(\eta)\Big
|\frac{ds}{s}\nonumber
\\&\le&\int_0^1 d|\mu|(\eta)\int_\eta^1 \frac{ds}{s}+\int_1^{\infty}d|\mu|(\eta)\int_1^\eta
\frac{ds}{s}\nonumber
\\&=&\int_0^\infty |\log \eta|\, d|\mu|(\eta) <\infty.
\nonumber\eea Similarly we have  $$\int_0^\infty k(s)
ds=\int_0^\infty \log \frac{1}{\eta} \, d\mu (\eta)=c_\mu,$$ which
gives (\ref{ka}). Thus, to complete the proof, it remains to justify
application of  Fubini's theorem leading to (\ref{eq:3100}). To this
end, it suffices to show that the repeated integral $$
\int_\e^\infty \frac{dt}{t}\int_0^\infty |Q_{t\eta} f(x)|\,d|\mu|
(\eta)$$ is finite for almost all $x$ in $\bbr^n$. We write it as
$A(x)+B(x)$, where
$$
A(x)=\int_\e^\infty \frac{dt}{t}\int_0^{1/t}|Q_{t\eta} f(x)|\,d|\mu|
(\eta), \quad B(x)=\int_\e^\infty
\frac{dt}{t}\int_{1/t}^\infty|Q_{t\eta} f(x)|\,d|\mu| (\eta).$$
Since the least radial decreasing majorant of  $q$ is integrable
(see property (b) in Definition \ref {d1}), then $\sup_{t>0}|Q_{t}
f(x)|\le c\, M_f (x)$ where $M_f (x)$ is the Hardy-Littlewood
maximal function, which is finite for almost $x$; see e.g.,
\cite[Theorem 2, Section 2, Chapter III]{St}. Hence, for almost $x$,
$$ A(x)\le c\,M_f (x)\int_\e^\infty \frac{dt}{t}\int_0^{1/t}d|\mu|
(\eta)= c\,M_f (x)\int_0^{1/\e}\big(\log \frac{1}{\eta}-\log
\e\big) d|\mu| (\eta)<\infty.$$ To estimate $B(x)$, we observe that
since $q\in L_r, \; r>1$, then, by Young's inequality
$$
||Q_t f||_s \le ||f||_p\, ||q_t||_r=t^{-\del}||f||_p\, ||q||_r,
\qquad \del=n(1-1/r)>0, \quad
\frac{1}{s}=\frac{1}{r}+\frac{1}{p}-1.$$ This gives $$\Big
\|\int_{1/t}^\infty|Q_{t\eta} f(x)|\,d|\mu| (\eta)\Big \|_s \le
t^{-\del}||f||_p\,\int_{1/t}^\infty \eta^{-\del}\,d|\mu| (\eta),$$
and therefore, \bea ||B||_s&\le&||f||_p\,\int_\e^\infty
\frac{dt}{t^{1+\del}}\int_{1/t}^\infty\eta^{-\del}\,d|\mu|
(\eta)\nonumber\\&=&\frac{||f||_p}{\del}\Big (\,
\int_0^{1/\e}d|\mu|(\eta)+\frac{1}{\e^\del}\int_{1/\e}^\infty
\eta^{-\delta} d|\mu|(\eta)\Big )\le \frac{||f||_p\,
||\mu||}{\del}<\infty. \nonumber \eea This completes the proof.
\end{proof}

\section{Wavelet Transforms Associated to One-parametric Semigroups and Inversion of Potentials}

In this section we consider an important subclass of wavelet
 transforms, generated by certain one-parametric semigroups of
 operators. Some composite wavelet
 transforms from the previous section belong to this subclass.

\subsection{Basic examples}
\begin{example}\label{e1} Consider the {\it Poisson semigroup} $\mathcal{P}_{t}$
generated by the Poisson integral
\begin{equation}
\mathcal{P}_{t}f(x)=\int_{\mathbb{R}^{n}}p(y,t)f(x-y)\,dy\text{ \ , \ }%
t>0  \label{1.1}
\end{equation}
 with the Poisson kernel
\begin{equation}
p(y,t)=\frac{\Gamma \left( (n+1)/2\right) }{%
\pi ^{(n+1)/2}}\frac{t}{(t^{2}+\left| y\right|
^{2})^{(n+1)/2}}=t^{-n}p(y/t, 1); \label{1.2}
\end{equation}
see \cite{31}, \cite {St}. In this specific case,  the kernel
function of the relevant composite wavelet transform is  $q(x)\equiv
p(x, 1)$ and on the Fourier transform side we have \be\label{puf}
F[\mathcal{P}_{t}f](\xi)=e^{-t\left| \xi\right| } Ff(\xi).\ee
\end{example}
\begin{example} Another important example is the {\it Gauss-Weierstrass
semigroup} $\mathcal{W}_{t}$ defined by
\begin{equation}
\mathcal{W}_{t}f(x)=\int_{\mathbb{R}^{n}}w(y,t)f(x-y)\,dy,\qquad
F[w(\cdot ,t)](\xi )=e^{-t|\xi |^{2}}, \quad t>0;  \label{2.4}
\end{equation}
see \cite{31}. The  Gauss-Weierstrass kernel $w(y,t)$ is explicitly
computed as \be w(y,t)=(4\pi t)^{-n/2}\exp (-\left| y\right|
^{2}/4t).\ee In comparison with (\ref{qu}), here the scaling
parameter $t$ is replaced by $\sqrt{t}$, so that
\be\label{scp2}w(y,t)=(\sqrt{t})^{-n}q (y/\sqrt{t}), \qquad q (y)=w
(y, 1)=(4\pi)^{-n/2}\exp (-\left| y\right| ^{2}/4), \ee and the
corresponding wavelet transform has the form \be \label{GWtr}
Wf(x,t)=\int_0^\infty \mathcal{W}_{t\eta} f(x)\, e^{-at\eta}\,
d\mu(\eta), \qquad  x \in \bbr^n, \;\;  t>0,\; \; a\ge 0.\ee This
 agrees with (\ref{cwtrw}) up to an obvious change of scaling
parameters.
\end{example}
\begin{example} The following interesting example does not fall into the scope
of wavelet transforms in Section 2, however, it has a very close
nature. Consider the {\it metaharmonic  semigroup} $\M_t$ defined by
\be\label{sg3}(\M_t f)(x) \! = \! \int\limits_{\mathbb R^n}
m(y,t)f(x \! - \! y)\,dy, \quad F[m(\cdot,t)](\xi) \! = \!
e^{-t\sqrt{1+|\xi|^2}};\ee see \cite[p. 257-258]{21}. The
corresponding kernel has  the form \be m
(y,t)=\frac{2t}{(2\pi)^{(n+1)/2}}\, \frac{K_{(n+1)/2} (\sqrt{|y|^2
+t^2})}{(\sqrt{|y|^2 +t^2})^{(n+1)/2}},\ee where
$K_{(n+1)/2}(\cdot)$ is the McDonald function. The
 relevant wavelet transform is  \be \label{mh}
Wf(x,t)=\int_0^\infty \mathcal{\M}_{t\eta} f(x)\, d\mu(\eta), \qquad
x \in \bbr^n, \;\;  t>0.\ee
\end{example}

This list of examples can be continued \cite{8}.

\subsection{Operators of the potential type}
One of the most remarkable applications of wavelet transforms
associated to  the Poisson, Gauss-Weierstrass, and metaharmonic
semigroups is that they pave the way to a series of explicit
inversion formulas for operators of the potential type arising in
analysis and mathematical physics. Typical examples of such
 operators are the following:
\begin{eqnarray}
\qquad I^{\alpha }f &=&F^{-1}\left| \xi \right| ^{-\alpha }Ff\equiv
(-\Delta
)^{-\alpha /2}f \quad \text{\rm (Riesz potentials)}, \label{1.8} \\
\qquad J^{\alpha }f &=&F^{-1}(1+\left| \xi \right| ^{2})^{-\alpha
/2}Ff\equiv
(E-\Delta )^{-\alpha /2}f \quad \text{\rm (Bessel potentials)}, \label{1.9} \\
\qquad \mathcal{F}^{\alpha }f &=&F^{-1}(1+\left| \xi \right|
)^{-\alpha }Ff\equiv (E+\sqrt{-\Delta })^{-\alpha }f  \quad
\text{\rm (Flett potentials)}.\label{1.10}
\end{eqnarray}
Here $Re\,\alpha >0$, \ $\left| \xi \right| =\left( \xi
_{1}^{2}+\cdots +\xi _{n}^{2}\right) ^{1/2},$ \ $\Delta
=\sum_{k=1}^{n}\frac{\partial ^{2}}{\partial x_{k}^{2}}$ is the
Laplacean, and $E$ is the identity operator. For \ $f\in
L_{p}(\mathbb{R}^{n}),$ \ $1\leq p<\infty ,$ these potentials have
 remarkable integral representations via the Poisson and Gauss-Weierstrass semigroups, namely,
\bea \label{pot1}I^{\alpha }f(x)&=&\frac{1}{\Gamma (\alpha
)}\int_{0}^{\infty }t^{\alpha -1}\mathcal{P}_{t}f(x)\,dt,\qquad
0<Re\,\alpha <n/p\,;\\ \label{pot2}J^\alpha f(x) &=& \frac{1}{\Gamma
(\alpha/2 )} \int_0^\infty
t^{\alpha/2-1}e^{-t}\,\W_t f(x) \,dt,\qquad 0<Re\,\alpha <\infty ;\\
\label {fla}\mathcal{F}^{\alpha }f(x)&=&\frac{1}{\Gamma (\alpha
)}\int_{0}^{\infty }t^{\alpha
-1}e^{-t}\mathcal{P}_{t}f(x)\,dt,\,\qquad 0<Re\,\alpha <\infty ;\eea
see \cite{30},  \cite{21},  \cite{12}. Regarding Flett potentials,
see, in particular, \cite[p. 446-447]{12}, \cite[p. 541-542]{28},
\cite{9}. We also mention another interesting representation of the
Bessel potential, which is due to Lizorkin \cite{19} and employs the
metaharmonic semigroup, namely,
\begin{equation}
J^{\alpha }f(x)=\frac{1}{\Gamma (\alpha )}\int_{0}^{\infty }t^{\alpha
-1}\mathcal{M}_{t}f(x)\,dt,\qquad 0<Re\,\alpha <\infty .  \label{1.13}
\end{equation}

Equalities (\ref{pot1})-(\ref{1.13}) have the same nature as
classical Balakrishnan's formulas for fractional powers of operators
(see  \cite[p. 121]{28}).

Let us show how these equalities generate wavelet inversion formulas
for the corresponding potentials. The core of the method is  the
following statement which is a particular case of Lemma 1.3 from
\cite{23}.
\begin{lemma}\label{lB} Given a finite Borel measure $\mu $ on $[0,\infty )$
and a complex number $\alpha ,$ \ $ \a'=Re\,\alpha \geq 0$, let
\begin{equation}
\lambda _{\alpha }(s)=s^{-1} I_+^{\alpha +1}\mu  (s),  \label{1.15}
\end{equation}
where
\begin{equation}
 I_+^{\alpha +1}\mu(s)=\frac{1}{\Gamma (\alpha +1)}%
\int_{0}^{s}(s-\eta )^{\alpha }d\mu (\eta )  \label{1.16}
\end{equation}
is the Riemann-Liouville fractional integral of order $\alpha +1$ of
the measure $\mu$. Suppose that $\mu $ satisfies the following
conditions: \bea \label{con1}&{}&\text{ \ }\int_{1}^{\infty
}\eta^{\gamma }d|\mu |(\eta
)<\infty \text{ \ \ \textit{for some} \ }\gamma >\a' ; \\
\label{con2}&{}&\text{ \ }\int_{0}^{\infty }\eta^{j}d\mu
(\eta)=0\text{ \  \ }\forall j=0,1,\ldots ,[Re\,\alpha ]\text{ \
\textit{(the integer part of} }\a' \text{).} \eea Then

\begin{eqnarray}\label{impo}
\lambda _{\alpha }(s)=\left\{
\begin{array}{lcl}
  O(s^{\a'-1}),   & \text{if} &    0<s<1,\\
O(s ^{-1-\del})\; \text{for some $\del >0$}, \; & \text{if} & s>1,
\end{array}
\right.
\end{eqnarray}
 and
\begin{eqnarray}
c_{\alpha,\mu } &\equiv &\int_{0}^{\infty }\lambda _{\alpha
}(s)\,ds=\int_{0}^{\infty }\frac{\tilde{\mu}(t)}{t^{\alpha +1}}\,dt
\notag \\
{}\nonumber\\
&=&\label{1.17n}\left\{
\begin{array}{lcl}\displaystyle{
\Gamma (-\alpha )\int_{0}^{\infty }\eta^{\alpha }\,d\mu (\eta)} &
\text{if} & \alpha \notin \mathbb{N}_{0}=\{0,1,2,\ldots \}, \\
{}\\
\displaystyle{\frac{(-1)^{\alpha +1}}{\alpha !}\int_{0}^{\infty
}\eta^{\alpha }\log \eta \, d\mu (\eta)} & \text{if} & \alpha \in
\mathbb{N}_{0},
\end{array}
\right.
\end{eqnarray}
where  $\tilde{\mu}(t)=\int_{0}^{\infty }e^{-t\eta}d\mu (\eta)$ is
the Laplace transform of $\mu$.
\end{lemma}

The estimate (\ref{impo}) is important in proving almost everywhere
convergence in forthcoming inversion formulas.

Consider, for example, Flett potential (\ref{1.10}), (\ref{fla}),
and make use of the composite wavelet transform \be \label {pwt}
W\vp (x,t)=\int_0^\infty
\mathcal{P}_{t\eta}\vp(x)\,e^{-t\eta}\,d\mu (\eta),\ee cf. Example
\ref{e1} and (\ref{cwtrw}) with $a=1$.
\begin{theorem}
\label{t1.5}\ Let $f\in L_{p},$ \ $1\leq p\leq \infty ,$ \
 and let $\varphi =\mathcal{F}^{\alpha
}f,$ \ $\alpha >0$, be the Flett potentials of $f.$ Suppose that
$\mu $ is a finite Borel measure on $[0,\infty )$ satisfying
(\ref{con1}) and (\ref{con2}).
 Then
\begin{equation}
\int_{0}^{\infty }W_{\mu }\varphi (x,t )\,\frac{dt }{t
^{1+\alpha }}\equiv \lim_{\varepsilon \rightarrow 0}\int_{\varepsilon
}^{\infty }W_{\mu }\varphi (x,t )\,\frac{dt }{t^{1+\alpha }}%
=c_{\alpha ,\mu }f(x),  \label{1.20}
\end{equation}
where \ $c_{\alpha ,\mu }$ is defined by (\ref{1.17n}) and the
limit is interpreted in the $L_{p}-$norm\ and pointwise a.e. on \
$\bbr^n$. If $f\in C_{0}$, the statement remains true with the
limit in (\ref {1.20}) interpreted in the sup-norm.
\end{theorem}
\begin{proof} We sketch the proof and address the reader to \cite{9} for
details. Changing the order of integration, owing to (\ref{pwt}),
(\ref{fla}), and the semigroup property of the Poisson integral, we
get \be W\varphi (x,t )=\frac{1}{\Gamma (\alpha )}\int_{0}^{\infty
}d\mu (\eta)\int_{t\eta}^{\infty }(\rho -t\eta )_{+}^{\alpha
-1}\,e^{-\rho }\,\mathcal{ P}_{\rho }f(x)\,d\rho . \label{1.21}\ee
Then further  calculations give \be \label {kuku}\int_{\varepsilon
}^{\infty }W\varphi (x,t )\frac{dt }{t^{1+\alpha }}=
\int_{0}^{\infty }e^{-\varepsilon s }\mathcal{P}_{\varepsilon s
}f(x)\,\lambda _{\alpha}(s) \,ds, \quad \lambda _{\alpha}\left(s
\right)= s^{-1}I_+^{\a +1} \mu (s),\ee cf. (\ref{1.15}). It remains
to applied Lemma \ref{lB} combined with the standard machinery of
approximation to the identity.
\end{proof}

Potentials (\ref{1.8})-(\ref{1.10}) and many others can be similarly
inverted  by making use of the  wavelet transforms associated with
suitable semigroups; see \cite{8}, \cite{9}.

\subsection{Examples of wavelet measures}
 Examples of  wavelet measures, that obey
the conditions of Lemma \ref{lB} with $c_{\alpha ,\mu }\neq 0$, are
the following.

1. Fix an integer $m>Re \,\a $ and choose an even Schwartz function
$h(\eta )$ \ on $\mathbb{R}^{1}$  so that $$h^{\left( k\right)
}(0)=0 \quad \forall  \,k=0,1,2,...,\quad \text{\rm and}\quad
\int_{0}^{\infty }\eta^{\a -m}h\left( \eta \right) d\eta \neq 0.$$
One can take, for instance, $h\left( \eta \right) =\exp \left( -\eta
^{2}-1/\eta ^{2}\right) ,$ \ $h\left( 0\right) =0.$ Set $d\mu \left(
\eta \right) =h^{\left( m\right) }\left( \eta \right) d\eta .$ It is
not difficult to show that $\int_{0}^{\infty }\eta ^{k}d\mu \left(
\eta \right) =0$, $\forall \ k=0,1,...,[Re \,\a]$, \ and $c_{\a,\mu
}\neq 0.$

2. Let $\mu =\sum\limits_{j=0}^{m}\binom{m}{j}\left( -1\right)
^{j}\delta _{j}$,  where $m>Re \,\a$ \ is a fixed integer and $
\delta _{j}=\delta _{j}\left( \eta \right) $ denotes the unit mass
at the point $\eta =j$, i.e., $\left\langle \delta
_{j},f\right\rangle =f(j).$ It is  known \cite[p. 117]{28}, \ that
\begin{equation*}
\int_{0}^{\infty }\eta ^{k}d\mu \left( \eta \right) \equiv
\sum\limits_{j=0}^{m}\binom{m}{j}\left( -1\right) ^{j}j^{k}=0,\text{ \ }%
\forall \text{ \ }k=0,1,...,m-1\text{ \ }\mathrm{(we\ set\ 0^{0}=1\ ).}
\end{equation*}
Moreover, $c_{\a ,\mu }=$ $\int_{0}^{\infty }t^{-\a -1}\left(
1-e^{-t}\right) ^{m}dt\neq 0.$

\section{Wavelet transforms with the generalized  translation operator}

Continuous wavelet transforms, studied in the previous sections,
 rely on the classical Fourier analysis on $\bbr^n$.
Interesting modifications of these  transforms and the corresponding
potential operators arise in the framework of the Fourier-Bessel
harmonic analysis associated to the Laplace-Bessel differential
operator \be \label{fb} \Delta_\nu =\sum\limits_{k=1}^{n}
\frac{\partial^2}{\partial
x_k^2}+\frac{2\nu}{x_n}\frac{\partial}{\partial x_n}\ , \qquad
\nu>0.\ee
 This analysis
amounts to pioneering works by Delsarte \cite{De} and Levitan
\cite{18}, and was extensively developed in subsequent publications;
see \cite{17}, \cite{32}, \cite{bh}, and references therein.

Let $\mathbb R^n_+=\{x: \, x=(x_1,\ldots , x_n)\in \mathbb R^n, \
x_n > 0\}$ and $x'=(x_1,\ldots , x_{n-1})$. Denote
$$
L_{p,\nu}(\bbr^n_+)=\Big \{f: ||f||_{p, \nu}=\Big (
\int_{\bbr^n_+}\left| f(x)\right| ^{p}x_n^{2\nu}dx\Big
)^{1/p}<\infty \Big \}.
$$
The Fourier-Bessel harmonic analysis is adopted to {\it the
generalized convolutions} \be\label{conv} (f \ast
g)(x)=\int_{\mathbb
 R^n_+} f(y)(T^y g)(x) \, y_n^{2\nu}dy, \qquad x \in \mathbb
 R^n_+ ,\ee
with the {\it generalized translation  operator} \be\label{gtop}
(T^yf)(x)=\frac{\Gamma(\nu +1/2)}{\Gamma(\nu)\Gamma(1/2)}\int_0^\pi
\!\! f(x'-y',\sqrt{x_n^2-2x_ny_n\cos \alpha
+y_n^2})\,\sin^{2\nu-1}\alpha \,d\alpha,  \ee
 \cite
{17}, \cite {18}, \cite {32}. The Fourier-Bessel transform $F_\nu$,
for which $F_{\nu }\left( f* g\right) =F_{\nu }\left( f\right)
F_{\nu }(g),$ \ is defined by
\begin{equation}\label{eq:211}
(F_\nu f)(\xi)=\int_{\mathbb R^n_+} f (x)\,e^{i\xi' \cdot
x'}j_{\nu-1/2}(\xi_n x_n) \, x_n^{2\nu} dx, \qquad \xi \in \mathbb
 R^n_+.
\end{equation}
Here $j_\lam (\tau)=2^\lam\Gamma (\lam+1) \, \tau^{-\lam}
J_\lam(\tau)$, where $J_\lam(\tau)$ is the Bessel function of the
first kind.  The {\it generalized Gauss-Weierstrass, Poisson, and
metaharmonic semigroups} $ \; \{\W_t^{(\nu)}\}, \; \{\P_t^{(\nu)}\},
\; \{\M_t^{(\nu)}\}$ are defined  as follows: \bea
 \label{k11}(\W_t^{(\nu)}f)(x)&=&\int\limits_{\mathbb
 R^n_+}w^{(\nu)}(y,t)(T^yf)(x) \, y_n^{2\nu}dy,\\
 && \hskip -1.5truecm F_\nu [w^{(\nu)}(\cdot,t)](\xi) =e^{-t|\xi|^2}; \nonumber \\
 \label{k21}(\P_t^{(\nu)}f)(x)&=&\int\limits_{\mathbb
 R^n_+}p^{(\nu)}(y,t)(T^yf)(x) \, y_n^{2\nu}dy, \\
  && \hskip -1.5truecm F_\nu [p^{(\nu)}(\cdot,t)](\xi) =e^{-t|\xi|}; \nonumber  \\
 \label{k31}  (\M_t^{(\nu)}f)(x)&=& \int\limits_{\mathbb
   R^n_+}m^{(\nu)}(y,t)(T^yf)(x) \, y_n^{2\nu}dy, \\
  &&  \hskip -1.5truecm F_\nu [m^{(\nu)}(\cdot,t)](\xi) =e^{-t\sqrt{1+|\xi|^2}}. \nonumber
\eea  The corresponding kernels $ w^{(\nu)}(y,t), \;
p^{(\nu)}(y,t)$, and $m^{(\nu)}(y,t)$ have the form \bea
\label{nu1}  w^{(\nu)}(y,t)&=& \frac{2\pi^{\nu +1/2}}{\Gamma(\nu +1/2)} \, (4\pi t)^{-(n+2\nu)/2} e^{-|y|^2/4t},\\
\label{nu2} p^{(\nu)}(y,t)&=& \frac{2 \,
\Gamma\big((n+2\nu+1)/2\big)}{\pi^{n/2}\Gamma (\nu+1/2)}\,
\frac{t}{(|y|^2
+t^2)^{(n+2\nu +1)/2}},\\
\qquad \label{nu3} m^{(\nu)}(y,t)&=& \frac{ 2^{-\nu+3/2}t}{\Gamma
(\nu+1/2)(2\pi)^{n/2}}\, \frac{K_{(n+2\nu +1)/2} (\sqrt{|y|^2
+t^2})}{(\sqrt{|y|^2 +t^2})^{(n+2\nu +1)/2}}.\eea
 More information
about these semigroups and their modifications
$$\{e^{-t}\W^{(\nu)}_t\}, \qquad   \{e^{-t}\P^{(\nu)}_t\}, \qquad \{
e^{-t}\M^{(\nu)}_t\},$$  can be found in \cite {2}, \cite {3}, \cite {13}.

Modified Riesz, Bessel, and Flett  potentials with the generalized
translation operator (\ref{gtop}) are formally defined in terms of
the Fourier-Bessel transform by
\begin{eqnarray}
I_{\nu }^{\alpha }f &=&F_{\nu }^{-1}\left| \xi \right| ^{-\alpha }F_{\nu
}f\equiv \left( -\Delta _{\nu }\right) ^{-\alpha /2}f,  \label{2.15} \\
\mathcal{J}_{\nu }^{\alpha }f &=&F_{\nu }^{-1}(1+\left| \xi \right|
^{2})^{-\alpha /2}F_{\nu }f\equiv \left( E-\Delta _{\nu }\right) ^{-\alpha
/2}f,  \label{2.17}\\
\mathcal{F}_{\nu }^{\alpha }f &=&F_{\nu }^{-1}(1+\left| \xi \right|
)^{-\alpha }F_{\nu }f\equiv \left( E+\sqrt{-\Delta _{\nu }}\right) ^{-\alpha
}f,  \label{2.16}
\end{eqnarray}
respectively. Here $Re\,\alpha >0$ and $\Delta _{\nu }$ is the
Laplace-Bessel differential operator (\ref{fb}). These generalized
potentials have analogous to (\ref{pot1})-(\ref{fla})
representations in terms of the semigroups
(\ref{k11})-(\ref{k31}), namely, if  $ f \in L_{p,\nu}(\bbr^n_+)$
then \bea \label{pot1n}I_{\nu }^{\alpha }f(x)&=&\frac{1}{\Gamma
(\alpha )}\int_{0}^{\infty }t^{\alpha
-1}\mathcal{P}_{t}^{(\nu)}f(x)\,dt,\qquad 0<Re\,\alpha
<(n+2\nu)/p\,,\\
J_{\nu }^\alpha f(x) &=& \frac{1}{\Gamma (\alpha/2 )} \int_0^\infty
t^{\alpha/2-1}e^{-t}\,\W_t^{(\nu)} f(x) \,dt,\qquad 0<Re\,\alpha <\infty ,\\
\label {flan}\mathcal{F}_{\nu }^{\alpha }f(x)&=&\frac{1}{\Gamma
(\alpha )}\int_{0}^{\infty }t^{\alpha
-1}e^{-t}\mathcal{P}_{t}^{(\nu)}f(x)\,dt,\,\qquad 0<Re\,\alpha
<\infty .\eea Moreover,
\begin{equation}
J_{\nu }^{\alpha }f(x)=\frac{1}{\Gamma (\alpha )}\int_{0}^{\infty }t^{\alpha
-1}\mathcal{M}_{t}^{(\nu)}f(x)\,dt,\qquad 0<Re\,\alpha <\infty .  \label{1.13n}
\end{equation}

We denote by $S_{t}^{(\nu)}$  any of the semigroups \be\label{sgr}
\W_t^{(\nu)}, \; \;e^{-t}\W_t^{(\nu)}, \;\;
\mathcal{P}_{t}^{(\nu)},\; \;e^{-t}\P^{(\nu )}_t, \;\;
\mathcal{M}_{t}^{(\nu)}, \;\; e^{-t}\M^{(\nu)}_t,\ee and define
the relevant wavelet transform (cf. (\ref{cwtr}))
\begin{equation}
\mathfrak{S}^{(\nu)}f(x,t)=\int_{0}^{\infty }S_{t\eta }^{(\nu)}\,f(x)\,
d\mu (\eta ), \qquad t>0,  \label{2.10}
\end{equation}
generated by a finite Borel measure $\mu$ on $[0,\infty )$.

There exist analogs of  Calder\'{o}n's reproducing formula for
wavelet transforms (\ref{2.10}) of functions belonging to the
weighted space $L_{p,\nu}(\bbr^n_+)$ and inversion formulas for
potentials $I_{\nu }^{\alpha }f, \; J_{\nu }^{\alpha }f,  \;
\mathcal{F}_{\nu }^{\alpha }f$, when $f\in L_{p,\nu}(\bbr^n_+)$.
For example, the following statement holds.
\begin{theorem}
\label{t2.6} Let $\varphi =I_{\nu }^{\alpha }f$, $f\in
L_{p,\nu}(\bbr^n_+)$, $1\leq p<\left( n+2\nu \right) /\alpha $, and
suppose that $\mu $ is a finite Borel measure on $[0,\infty )$
satisfying (\ref{con1}) and (\ref{con2}). If \
$\mathfrak{S}^{(\nu)}\varphi $ is the wavelet transform of $\varphi$
associated with the generalized Poisson semigroup
$\mathcal{P}_{t}^{(\nu)}$, then
\begin{equation}
\int_{0}^{\infty }\frac{\mathfrak{S}^{(\nu)}\varphi (x,t)}{t^{1+\alpha
}}dt=\lim_{\varepsilon \rightarrow 0}\int_{\varepsilon }^{\infty }%
\frac{\mathfrak{S}^{(\nu)}\varphi \left( x,t\right) }{t^{1+\alpha }}%
dt=c_{\alpha ,\mu }f(x),  \label{2.19}
\end{equation}
where $c_{\alpha ,\mu }$ is defined by (\ref{1.17n}). The limit in
(\ref{2.19}) exists in the $L_{p,\nu}(\bbr^n_+)$-norm  and in the
a.e. sense. If $f\in C_{0}$, the convergence in (\ref{2.19}) is
uniform.
\end{theorem}
The proof of this theorem is presented in \cite{8} in the general
context of the so-called admissible semigroups. This context
includes all semigroups (\ref{sgr}).

\section{ Beta-semigroups}

We remind basic formulas from Section 3.1 for the kernels of the
Poisson and Gauss-Weierstrass semigroups: \be\label{beta1}
F[p(\cdot ,t)](\xi )=e^{-t\left| \xi\right| }, \qquad p(y,t)=\frac{\Gamma \left( (n+1)/2\right) }{%
\pi ^{(n+1)/2}}\frac{t}{(t^{2}+\left| y\right| ^{2})^{(n+1)/2}};\ee
\be \label{beta2} F[w(\cdot ,t)](\xi )=e^{-t|\xi |^{2}}, \qquad
w(y,t)=(4\pi t)^{-n/2}\exp (-\left| y\right| ^{2}/4t).\ee It would
be natural to consider a more general semigroup generated by the
kernel $w^{(\b)}(y,t)$ defined by \be \label{beta} F[w^{(\b)}(\cdot
,t)](\xi )=e^{-t|\xi |^{\b}}, \qquad \b>0.\ee This semigroup arises
in diverse contexts of analysis, integral geometry,  and
probability; see, e.g., \cite{Fe}, \cite{Ko}, \cite{La}, \cite{R8}.
Unlike (\ref{beta1}) and (\ref{beta2}), the kernel function
$w^{(\b)}(y,t)$ cannot be computed explicitly, however, by taking
into account that \be w^{(\b)}(y,t) =t^{-n/\beta } w^{(\b)}(
t^{-1/\beta }y),\qquad w^{(\b)}(y)\equiv w^{(\b)}(y,1), \ee
 properties of $w^{(\b)}(y,t)$ are well
  determined by the following lemma.
\begin{lemma}\label{lb}  The function
 \be\label{gaql}
w^{(\b)}(y) =F^{-1}[e^{-|\cdot|^\b}](y)=(2\pi)^{-n}
\int_{\rn}e^{-|\xi|^\b} e^{i y\cdot \xi}\, d\xi, \qquad \b>0,\ee is
uniformly continuous on $\rn$. If $\b$ is an even integer, then
 $w^{(\b)}(y)$ is infinitely smooth
 and rapidly decreasing. More generally, if $\b\neq 2,4,\ldots $, then $w^{(\b)}(y)$
 has the following behavior when $|y| \to \infty$: \be\label{cbe}
w^{(\b)}(y) =c_\b |y|^{-n-\b} (1+o(|y|)), \quad
c_\b=-\frac{2^{\b}\pi^{-n/2} \Gam ((n+\b)/2)}{ \Gam (-\b/2)}.\ee
If $0<\b\le 2$, then $w^{(\b)}(y)>0$ for all $y \in \rn$.
\end{lemma}
\begin{proof} (Cf. \cite[p. 44, for $n=1$]{Ko}). The uniform
continuity of $w^{(\b)}(y)$ follows immediately from (\ref{gaql}).
Note that if $\b$ is an even integer, then $e^{-|\cdot|^\b}$ is a
Schwartz
 function and therefore, $w^{(\b)}(y)$ is infinitely smooth
 and rapidly decreasing.
Let us prove positivity of $w^{(\b)}(y)$ when $0<\b\le 2$.  For
$y=0$ and for the cases $\b=1$ and  $\b=2$, this is obvious.
 Let
$0<\b< 2$. By Bernstein's theorem \cite[Chapter 18, Sec. 4]{Fel},
there is a non-negative finite measure $\mu_\b$ on $[0,\infty)$ so
that $ e^{-z^{\b/2}}=\int_0^\infty e^{-tz}\,d\mu_\b (t)$, $z\in
[0,\infty)$. Replace $z$ by $|\xi|^2$ to get \be\label{751}
e^{-|\xi|^{\b}}=\int_0^\infty e^{-t|\xi|^2}\,d\mu_\b (t).\ee Then
the equality \be\label{75}
[e^{-t|\cdot\,|^2}]^{\wedge}(y)=\pi^{n/2}t^{-n/2}e^{-|y|^2/4t},
\qquad t>0,\ee yields \bea
 (2\pi)^{n}\,w^{(\b)}(y)&=&\int_{\rn}e^{i \xi\cdot
 y}d\xi\int_0^\infty e^{-t|\xi|^2}\,d\mu_\b (t)= \int_0^\infty d\mu_\b
(t)\int_{\rn}e^{i \xi\cdot y} e^{-t|\xi|^2}\,d\xi\nonumber\\&=&
\pi^{n/2}\int_0^\infty t^{-n/2}e^{-|y|^2/4t}\,d\mu_\b
(t)>0.\nonumber\eea The Fubini theorem is applicable here, because,
by (\ref{751}),
$$
\int_{\rn}|e^{i \xi\cdot y}|d\xi\int_0^\infty
e^{-t|\xi|^2}\,d\mu_\b  (t)=\int_{\rn}
e^{-|\xi|^{\b}}d\xi<\infty.$$

Let us prove (\ref{cbe}).  It suffices to show
 that
 \be\label{73p}
\lim\limits_{|y| \to \infty}|y|^{n
+\b}w^{(\b)}(y)=2^{\b}\pi^{-n/2-1}\Gam (1+\b/2)\Gam ((n+\b)/2)\,
\sin (\pi \b/2)\ee (we leave to the reader to check that the
right-hand side coincides with $c_\b$).  For $n=1$, this statement
can be found in \cite [Chapter 3, Problem 154]{PS} and  in \cite
[p. 45]{Ko}. In the general case, the proof is more sophisticated
and relies on the properties of Bessel functions. By the
well-known formula for the Fourier transform of a radial function
(see, e.g., \cite{31}), we write
$(2\pi)^{n}\,w^{(\b)}(y)=I(|\eta|)$, where \bea
I(s)&=&(2\pi)^{n/2}s^{1-n/2}\int_0^\infty e^{-r^{\b}} r^{n/2}
J_{{n/2-1}} (rs)\, dr\nonumber\\&=&(2\pi)^{n/2}s^{-n}\int_0^\infty
e^{-r^{\b}} \frac{d}{dr}\,[(rs)^{n/2} J_{{n/2}} (rs)]\,
dr.\nonumber\eea Integration by parts yields $$
I(s)=\b(2\pi)^{n/2}s^{-n/2}\int_0^\infty e^{-r^{\b}} r^{n/2+\b-1}
J_{{n/2}} (rs)\, dr.$$ Changing variable $z=s^\b r^\b$, we obtain
$$ s^{n +\b}I(s)=(2\pi)^{n/2} A(s^{-\b}), \qquad A(\del)=
\int_0^\infty e^{-z\del} z^{n/2\b} J_{n/2} (z^{1/\b})\,dz.
$$
We actually have to compute the limit $A_0=\lim\limits_{\del \to 0}
A(\del)$. To this end, we invoke Hankel functions $H_\nu^{(1)} (z)$,
so that $ J_\nu (z)=Re \,H_\nu^{(1)} (z)$ if $z$ is real \cite{Er}.
 Let $h_\nu (z)=z^\nu H_\nu^{(1)} (z)$. This is a single-valued
analytic function in the $z$-plane with cut $(-\infty, 0]$. Using
the  properties of the Bessel functions \cite{Er}, we  get
 \be\label {as}
\lim\limits_{z \to 0}h_\nu (z)=2^\nu \Gam (\nu)/\pi i,\ee \be\label
{as1} h_\nu (z) \sim \sqrt{2/\pi} \, z^{\nu -1/2}e^{iz-\frac{\pi i
}{2}(\nu +\frac{1}{2})}, \qquad z \to \infty.\ee Then we write
$A(\del)$ as $ A(\del)= Re \,\int_0^\infty e^{-z\del} h_{n/2}
(z^{1/\b})\,dz$ and change the line of integration from $[0,\infty)$
to $n_\theta=\{z: z=re^{i\theta}, \; r>0\}$ for small $\theta<\pi
\b/2$. By Cauchy's theorem, owing to (\ref{as}) and (\ref{as1}), we
obtain $ A(\del)= Re \,\int_{n_\theta} e^{-z\del} h_{n/2}
(z^{1/\b})\,dz$. Since for $z=re^{i\theta}$, $ h_{n/2}
(z^{1/\b})=O(1)$ when $r=|z|\to 0$ and $ h_{n/2} (z^{1/\b})=O(r^{(n
-1)/2\b} e^{-r^{1/\b}\sin (\theta /\b)})$ as $r\to \infty$, by the
Lebesgue theorem on dominated convergence, we get $ A_0=Re
\,\int_{n_\theta}  h_{n/2} (z^{1/\b})\,dz$. To evaluate
 the last integral, we again use analyticity and replace
$n_\theta$ by $n_{\pi \b/2}=\{z: z=re^{i\pi \b/2}, \; r>0\}$ to get
$$ A_0=Re \,\Big [e^{i\pi \b/2}\int_0^\infty   h_{n/2}
(r^{1/\b}e^{i\pi/2})\,dr\Big ].$$ To finalize calculations, we
invoke McDonald's function $K_\nu (z)$ so that
$$
h_\nu (z)=z^\nu H_\nu^{(1)} (z)=-\frac{2i}{\pi}(z e^{-i\pi/2})^\nu
K_\nu (z e^{-i\pi/2}).$$ This gives
$$
A_0=\frac{2\b}{\pi}\, \sin (\pi \b/2) \int_0^\infty s^{n/2 +\b-1}
K_{n/2} (s)\, ds.
$$ The last integral can be explicitly evaluated by the formula 2.16.2
(2) from \cite {PBM}, and we obtain the result.
\end{proof}

The Beta-semigroup $\mathcal{B}_{t}$ generated by the kernel
$w^{(\b)} (y,t)$ (see (\ref{beta})) is defined by
\begin{equation}
\mathcal{B}_{t}f(x)=\int_{\mathbb{R}^{n}}w^{\left( \beta \right)
}\left( y,t\right) f\left( x-y\right) dy, \qquad  t>0, \label{2.21}
\end{equation}
and the corresponding weighted wavelet transform has the form
\be\label{bew} W_a f(x,t)=\int_0^\infty  \B_{t\eta} f (x)\,
e^{-at\eta}\,d\mu (\eta),\ee
 where $a\ge 0$ is a fixed number which is in our disposal; cf
 (\ref{cwtrw}).
 Following \cite{A}, we introduce Beta-potentials
 \be\label{bpot}
 J_\b^\a f = (E+ (-\Del^{\b/2}))^{-\a/\b}f,\qquad \alpha >0, \quad \b >0,\ee
 that can be realized through the Beta-semigroup as
 \be\label{bpott}
J_\b^\a f (x)=\frac{1}{\Gamma (\alpha /\beta )}
\int_{0}^{\infty }t^{\a/\b-1}\text{ }e^{-t}\text{ }%
\mathcal{B}_{t}f(x)\text{ }dt. \ee For $\b=2$, (\ref{bpot})
coincides with the classical Bessel potential (\ref{1.9}), and
(\ref{bpott}) mimics (\ref{pot2}). Similarly, for $\b=1$, the
Beta-potentials coincide with the Flett potential  (\ref{fla}).

Explicit inversion formulas for Beta-potentials  can be obtained
with the aid of the  wavelet transform (\ref{bew}) as follows.
\begin{theorem}
\label{t2.7} Let $f\in L_{p}(\mathbb{R}^{n})$, \ $1\leq p< \infty $,
$\alpha >0, \; \b>0$. Suppose that $\mu $ is a finite Borel measure
on $[0,\infty )$ satisfying
\begin{eqnarray*}
&\text{(a) \ }&\int_{1}^{\infty }\eta^{\gamma }\text{
}d\left| \mu \right| \left( \eta\right) <\infty \text{ \ \ for some
\ }\gamma
>\alpha /\beta
; \\
&\text{(b) \ }&\int_{0}^{\infty }\eta^{j}\text{ }d\mu \left( \eta\right) =0, \quad \forall \,j=0,1,...,[\alpha /\beta ].
\end{eqnarray*}
If $\varphi =J_{\beta }^{\alpha }f$, then \be\label{form}
\int_{0}^{\infty }W\varphi \left( x,t\right) \frac{dt}{t^{1+\alpha
/\beta }}\equiv \lim_{\varepsilon \rightarrow 0}\int_{\varepsilon
}^{\infty }W\varphi \left( x,t\right) \frac{dt}{t^{1+\alpha /\beta
}}=c_{\alpha /\beta ,\mu }f(x), \ee where $c_{\a/\b,\mu }$ \ is
defined by (\ref{1.17n}) (with $\a$ replaced by $\a/\b$). The limit
in (\ref{form}) exists in the $L_{p}$-norm and pointwise
 for almost all x. If $f\in C_{0}$,  the convergence is uniform.
\end{theorem}

The proof of this theorem  mimics that of Theorem \ref {t1.5}; see
\cite {A} for details.

\begin{remark} The classical Riesz potential $I^{\alpha}f$ has an
integral representation via the Beta-semigroup, namely,
\be\label{form} I^{\alpha}f(x)= \frac{1}{\Gamma (\alpha /\beta )}
\int_{0}^{\infty }t^{{\a/\b}-1}\text{ } \mathcal{B}_{t}f(x)\text{
}dt.\ee Here $f \in L_{p}(\mathbb{R}^{n}), \ 1\leq p < \infty$ , and
\ $ 0< Re \alpha < n/p$. For the  cases $\beta=1$ and $\beta=2$ we
have the representations in terms of the Poisson and
Gauss-Weierstrass semigroups, respectively.

The potential $I^{\alpha}f$ can be inverted in the framework of the
$L_{p}$-theory by making use of (\ref{form}) and the composite
wavelet transform (\ref{bew}) with $a=0$.

\end{remark}

\section{Parabolic Wavelet Transforms}

The following anisotropic wavelet transforms of the composite type,
associated with the heat operators \be\label{Ho} \partial /\partial
t -\Delta, \qquad
 E+\partial /\partial t  -\Delta, \ee
were introduced by  Aliev and  Rubin \cite{7}. These transforms are
constructed using the Gauss-Weierstrass kernel $w(y, t) = (4\pi
t)^{-n/2} \exp (- |y|^2/4t)$ as follows. Let $\bbr^{n+1}$ be the
$(n+1)$-dimensional Euclidean space of points $(x, t)$, $x = (x_1,
\ldots, x_n) \in \bbr^n, \, $ $t \in \bbr^1$. We pick up a wavelet
measure $\mu$ on $[0,\infty)$, a scaling parameter $a>0$, and set
\be \label{pwtr} P_\mu f(x, t;a) = \int_{\bbr^n \times (0, \infty)}
f(x-\sqrt{a} y, t-a\tau) \,w(y, \tau) \, dyd\mu (\tau), \ee \be
\label{pwtrw}\P_\mu f(x, t; a) = \int_{\bbr^n \times (0, \infty)}
f(x - \sqrt{a} y, t- a\tau) \,w(y, \tau)\, e^{-a\tau}
\,dyd\mu(\tau)\ee (to simplify the notation, without loss of
generality we can assume $\mu (\{0\})=0$). We call (\ref{pwtr}) and
(\ref{pwtrw}) the {\it parabolic wavelet transform} and  the {\it
weighted parabolic wavelet transform}, respectively.

Parabolic potentials $H^\a f$ and $\H^\a f$, associated to
differential operators in (\ref{Ho}), are defined in the Fourier
terms by \be\label{ppo} F[H^\a f](\xi, \tau) = (|\xi|^2 + i
\tau)^{-\a/2} F[f] (\xi, \tau),\ee \be\label{ppo1}  F[\Cal H^\a
f](\xi, \tau) = (1+ |\xi|^2 + i \tau)^{-\a/2} F[f](\xi, \tau),\ee
where $F$ stands for the Fourier transform in $\bbr^{n+1}$. These
potentials were introduced by  Jones \cite{Jo} and Sampson \cite{Sa}
and used  as a tool for characterization of anisotropic function
spaces of fractional smoothness; see \cite{7} and references
therein. For $\a>0$, potentials $H^\a f$ and $\H^\a f$ are
representable by the integrals \bea \qquad H^\a f(x, t) &=& {1\over
\Gamma(\a/2)} \int_{\bbr^n \times (0,\infty)} \tau^{\a/2-1} w(y,
\tau) \,f(x-y, t-\tau)\,  dyd\tau,\\
 \H^\a f(x, t)&=&{1\over \Gamma (\a/2)} \int_{\bbr^n \times
(0, \infty)}  \tau^{\a/2-1} e^{-\tau} w(y, \tau) \,f(x-y, t-\tau) \,
dyd\tau.\eea Their behavior on functions $f \in L_p \equiv L_p
(\bbr^{n+1})$ is characterized by the following theorem.

\begin{theorem} \cite {Ba}, \cite {Ra}
\newline
{\rm I.} \ Let $f \in L_{p}, \; 1 \le p < \infty, \; 0 < \a <
(n+2)/p, \quad q = (n+2-\a p)^{-1} (n+2) p$.

{\rm (a)} \ The integral $(H^\a f)(x, t)$ converges absolutely for
almost all $(x, t) \in \bbr^{n+1}$.

{\rm (b)} \ For $p > 1$, \ \ the operator $H^\a$ is bounded from
$L_{p}$ into $L_{q}$.

{\rm (c)} \  For $p = 1$, $H^\a$ is an operator of the weak $(1, q)$
type:
$$
|\{ (x, t): |(H^\a f) (x, t) | > \gamma\}| \le \left({c\| f \|_1
\over \gamma}\right)^q.
$$
\newline
{\rm II.} \ The operator $\H^\a$ is bounded on $L_{p}$ for all $\a
\ge 0, \quad 1 \le p \le \infty$.
\end{theorem}

Explicit inversion formulas for parabolic potentials in terms of
wavelet transforms (\ref{pwtr}) and (\ref{pwtrw})  are given by the
following theorem.
\begin{theorem}  \cite{7} Let $\mu$ be a finite Borel measure on $[0,
\infty)$ satisfying the following conditions:

\bea \label{conn1}&{}&\text{ \ }\int_{1}^{\infty }\t^{\gamma
}d|\mu |(t
)<\infty \text{ \ \ \textit{for some} \ }\gamma >\alpha/2 ; \\
\label{conn2}&{}&\text{ \ }\int_{0}^{\infty }t^{j}d\mu (t)=0\text{
\ , \ }\forall j=0,1,\ldots ,[\alpha/2 ]. \eea Suppose that $\vp =
H^\a f, \; \; f \in L_p, \; \; 1 \le p < \infty, \; \; 0 < \a <
(n+2)/p$. Then

\be\label {inf}\int_0^\infty P_\mu \vp(x, t; a) \; {da\over
a^{1+\a/2}} \, \equiv \, \lim_{\e \to 0} \int^\infty_\e (\dots) =
c_{{\a}/2, \mu} \ f(x, t),\ee where $c_{{\a}/2, \mu}$ is defined
by (\ref{1.17n}) (with $\alpha$ replaced by $\alpha /2$).

The limit in (\ref {inf}) is  interpreted in the $L_{p}$-norm for
$1 \le p < \infty$ and a.e. on $\bbr^{n+1}$ for $1 < p < \infty$.

The same statement holds for all $\a > 0$ and $1 \le p \le \infty$
($L_\infty$ is identified with $C_0$)
 provided that $H^\a$ and $P_\mu$ are
replaced by $\H^\a$ and $\P_\mu$, respectively.
\end{theorem}

More general  results for parabolic wavelet transforms  with the
generalized translation associated to singular heat operators
\be\label{Hos} \partial /\partial t -\Delta_{\nu}, \qquad
 E+\partial /\partial t  -\Delta_{\nu}, \qquad \qquad\Big (\Delta _{\nu }=\sum\limits_{k=1}^{n}\frac{\partial ^{2}}{\partial
x_{k}^{2}}+\frac{2\nu}{x_{n}}\,\frac{\partial}{\partial x_{n}}\Big
),\ee
  were obtained in \cite{6}. These include the Calder\'{o}n-type reproducing
formula and explicit $L_p$-inversion formulas for parabolic
potentials with the generalized translation defined by
\begin{eqnarray*}
H_{\nu }^{\alpha }f(x,t) &=&F_{\nu }^{-1}[(\left| x\right| ^{2}+it)^{-\alpha
/2}\text{ }F_{\nu }f(x,t)], \\
\mathcal{H}_{\nu }^{\alpha }f(x,t) &=&F_{\nu }^{-1}[(1+\left| x\right|
^{2}+it)^{-\alpha /2}\text{ }F_{\nu }f(x,t)].
\end{eqnarray*}
In the last two expressions, $x\in \mathbb{R}_{+}^{n}=\{x\in
\mathbb{R}^{n}:$ \ $ x_{n}>0\}$,  $\; t\in \mathbb{R}^{1}$, and
$F_{\nu}$ \ is the Fourier-Bessel transform, i.e., the Fourier
transform with respect to the variables $ t $ and $x^{\prime
}=(x_{1},...,x_{n-1}),$ and the Bessel transform with respect to
$x_{n}>0.$ These results were applied in \cite{6, 7} to wavelet-type
 characterization of the parabolic Lebesque spaces.

\section{Some Applications to Inversion of the $k$-plane Radon Transform}

We recall some basic definitions. More information can be found in
\cite{{GGG}, 15, 10, 24, 25}. Let $\ \mathcal{G}_{n,k}$ \ and $
G_{n,k}$ be the affine Grassmann manifold of all non-oriented
$k$-dimensional planes ($k$-planes) $\tau $ \ in $\mathbb{R}^{n}$
and the
ordinary Grassmann manifold of $k$-dimensional linear subspaces $\zeta $ of $%
\mathbb{R}^{n}$, respectively. Each $k$-plane $\tau \in
\mathcal{G}_{n,k}$ is parameterized as $\tau =\left( \zeta \text{,
}u\right) $, where $\zeta \in G_{n,k}$ and $u\in \zeta ^{\perp }$
(the orthogonal complement of $\zeta $ in $\mathbb{R}^{n}$). We
endow $\mathcal{G}_{n,k\text{ }}$  with the product measure $d\tau
=d\zeta du$,  where  $d\zeta $ is the $O(n)$-invariant measure on
$G_{n,k\text{ }}$ of total mass $1,$ and $du$ denotes the Euclidean
volume element on $\zeta ^{\perp }$. The \textit{\ k-plane Radon
transform }of a function $f$ on $\mathbb{R}^{n}$ is defined by
\begin{equation}
\hat f (\tau )\equiv \hat f (\zeta \text{, }%
u)=\int_{\zeta }f(y+u)\,dy,  \label{3.1}
\end{equation}
where $dy$ is the induced Lebesque measure on the subspace \ $\zeta
\in G_{n,k}.$ \  This transform assigns to a function $f$ a
collection of integrals of $f$  over all $k$-planes in $\rn$. The
corresponding \textit{dual $k$-plane transform} of a function
$\varphi$ on $ \mathcal{G}_{n,k}$ is defined as the mean value of
$\varphi \left( \tau \right) $ over all $k$-planes $\tau $ through
$x\in \mathbb{R}^{n}$:
\begin{equation}
\check{\varphi}\left( x\right) =\int_{O(n)}\varphi (\mathcal{\sigma }\zeta
_{0}+x) \, d\mathcal{\sigma }, \qquad x\in \mathbb{R}^{n}.  \label{3.2}
\end{equation}
Here $\zeta _{0}\in G_{n,k}$ \ is an arbitrary fixed $k-$plane
through the origin. If $f\in L_{p}(\mathbb{R}^{n}),$ then $\hat f $
is finite a.e. on $\mathcal{G}_{n,k}$ \ if and only if $1\leq
p<n/k.$

Several inversion procedures are known for $\hat f$. One of the most
popular, which amounts to Blaschke and Radon, relies on  the Fuglede
formula \cite[p. 29]{15},
\begin{equation}
( \hat f)^{\vee }=d_{k,n}I^{k}f, \qquad
d_{k,n}=\left( 2\pi \right) ^{k}\sigma _{n-k-1}/\sigma _{n-1},
\label{3.3}
\end{equation}
and reduces reconstruction of $f$  to
inversion of the Riesz potentials $I^{k}f.$ The latter can also be inverted in many number of ways \cite{29}, \cite%
{28}, \cite{21}. In view of considerations in Section 3.2 and 5, one
can employ a composite wavelet transform generated by the Poisson,
Gauss-Weierstrass, or Beta semigroup and thus obtain new inversion
formulas for the $k$-plane transform on $\rn$ in terms of a wavelet
measure on the one-dimensional set $[0,\infty)$. For instance, this
way leads to the following
\begin{theorem}
\label{t3.1} Let $\varphi =\hat f$ be the $k$-plane Radon transform
of a function
$f\in L_{p}$, $1\leq p<n/k$. $\ $Let $\mu $ be a finite Borel measure on $%
[0,\infty )$ satisfying
\begin{eqnarray*}
&\text{(a) \ }&\int_{1}^{\infty }\eta ^{\gamma }d\left| \mu \right|
\left( \eta \right) <\infty \text{ \ for some }\gamma >k; \\
&\text{(b) \ }&\int_{0}^{\infty }\eta ^{j}d\mu \left( \eta \right) =0
\quad \forall \text{ }j=0,1,...,k.
\end{eqnarray*}
Let $W\check{\varphi}$ be the wavelet transform of
$\check{\varphi}$,  associated with the Poisson semigroup
(\ref{1.1}), namely, \be \label{pWtr}
W\check{\varphi}(x,t)=\int_0^\infty \mathcal{P}_{t\eta}
\check{\varphi}(x)\,  d\mu(\eta), \qquad  x \in \bbr^n, \;\;
t>0.\ee
 Then
\begin{equation}
\int_{0}^{\infty }W\check{\varphi}\left( x,t\right)
\frac{dt}{t^{1+k}}\equiv \lim_{\varepsilon \rightarrow \infty
}\int_{\varepsilon }^{\infty }W\check{\varphi}%
\left( x,t\right) \frac{dt}{t^{1+k}}=c_{k,\mu }f(x),  \label{3.4}
\end{equation}
where (cf. (\ref{1.17n})),
\begin{equation*}
c_{k,\mu }=\frac{\left( -1\right) ^{k+1}}{k!}\int_{0}^{\infty
}t^{k}\log t\text{ }d\mu \left( t\right) .
\end{equation*}
The limit in (\ref{3.4}) exists in the $L_{p}$-norm and pointwise
almost everywhere. If $f\in C_{0}\cap L_{p}$, the convergence is
uniform on $\rn$.
\end{theorem}

\begin{remark} The following observation might be interesting. Let \be
I_-^\a u(t)=\frac{1}{\Gam (\a)}\int_t^\infty (t-s)^{\a -1} u(s)\,ds,
\qquad t>0,\ee be the Riemann-Liouville integral of $u$. It is known
\cite[formula (16.9)]{21} that the Poisson integral takes the Riesz
potential $I^\a f$ to the Riemann-Liouville integral of the function
$t \to \P_t f$, namely,  \be \label {comb4} \P_t I^\a f=I_-^\a
\P_{(\cdot)} f.\ee
 Denoting by $R$ and $R^{\ast }$ the Radon $k$-plane transform and its dual, owing to  Fuglede's formula
 (\ref{3.3}), we have
\begin{equation}
R^{\ast }Rf=d_{k,n}\,I^{k}f.  \label{3.5}
\end{equation}
Combining (\ref{3.5}) and (\ref{comb4}), we get \begin{equation}
R_{t}^{\ast }Rf=d_{k,n}\,I_-^k \P_{(\cdot)} f, \qquad R_{t}^{\ast
}\varphi (x)=(\mathcal{P}_{t}R^{\ast }\varphi )(x) . \label{3.7}
\end{equation}
This formula has the same nature as the following one in terms of
the  spherical means, that lies in the scope of the classical
Funk-Radon-Helgason theory: \be\label{frh} (\hat f)^\vee_r
(x)=\sig_{k-1}\int_r^\infty (\frM_t f)(x) (t^2 -r^2)^{k/2 -1} t \,
dt; \ee see Lemma 5.1 in \cite{24}. Here $\sig_{k-1}$ is the volume
of the $(k-1)$-dimensional unit sphere, \be (\frM_t
f)(x)=\frac{1}{\sig_{n-1}}\int_{S^{n-1}} f(x+t\theta ) \,d \theta,
\qquad t>0, \ee and $(\hat f)^\vee_r (x)$ is the so-called {\it
shifted dual  $k$-plane transform}, which is the mean value of $\hat
f (\t)$ over all $k$-planes $\t$ at distance $r$ from $x$.

\end{remark}

\section{Higher-rank Composite Wavelet Transforms and Open Problems}

Challenging perspectives and  open problems for composite wavelet
transforms  are connected with functions of matrix argument and
their application to integral geometry. This relatively new area
encompasses the so-called higher-rank problems, when traditional
scalar notions, like distance or scaling, become matrix-valued.

\subsection{Matrix spaces, preliminaries}
We  remind basic notions, following \cite{R9}. Let $\frM_{n,m}  \sim
\bbr^{nm}$ be the space of real matrices $x=(x_{i,j})$ having $n$
rows and $m$
 columns, $n\geq m$;  $dx=\prod^{n}_{i=1}\prod^{m}_{j=1}
 dx_{i,j}$ is the volume element on $\Ma$, $x'$ denotes the transpose of  $x$, and  $I_m$
   is the identity $m \times m$
  matrix. Given a square matrix $a$,  we denote by
  $\det(a)$
  the determinant of $a$, and by $|a|$ the absolute value of
  $\det(a)$;
  $\tr (a)$ stands for the trace of $a$. For $x\in \Ma$, we denote  \be\label{xm}|x|_m =\det (x'x)^{1/2}.\ee
If $m=1$,  this is the usual Euclidean norm on $\bbr^n$. For $m>1$,
$|x|_m$ is the volume of the parallelepiped spanned  by the
column-vectors of  $x$.
   We use  standard
 notations $O(n)$   and $SO(n)$ for the orthogonal group and the
 special orthogonal group of $\bbr^{n}$ with the  normalized
 invariant measure of total mass 1. Let $\S_m \sim \bbr^{m(m+1)/2}$  be the space of $m \times m$ real
symmetric matrices $s=(s_{i,j})$
 with the volume element $ds=\prod_{i \le j} ds_{i,j}$. We denote by  $\p$  the cone of
positive definite matrices in $\S_m$; $\cpm$  is the closure of
$\p$, that is, the set of  all positive semi-definite $m\times m$
matrices.  For $r\in\p$ ($r\in\cpm$), we write $r>0$ ($r\geq 0$).
Given
 $a$ and  $b$ in  $S_m$, the inequality $a >b$  means $a - b \in
 \p$ and  the
symbol $\int_a^b f(s) ds$ denotes
 the integral over the set $(a +\p)\cap (b -\p)$.

  The group $G=GL(m,\bbr)$  of
 real non-singular $m \times m$ matrices $g$ acts transitively on $\p$
  by the rule $r \to g rg'$.  The corresponding $G$-invariant
 measure on $\p$ is  \be\label{2.1}
  d_{*} r = |r|^{-d} dr, \qquad |r|=\det (r), \qquad d= (m+1)/2 \ee
  \cite[p. 18]{Te}.

\begin{lemma}\label{12.2} \cite[pp. 57--59] {Mu}\hskip10truecm

\noindent
 {\rm (i)} \ If $ \; x=ayb$ where $y\in\Ma, \; a\in  GL(n,\bbr)$, and $ b \in  GL(m,\bbr)$, then
 $dx=|a|^m |b|^ndy$. \\
 {\rm (ii)} \ If $ \; r=q'sq$ where $s\in S_m$, and $q\in  GL(m,\bbr)$,
  then $dr=|q|^{m+1}ds$. \\
  {\rm (iii)} \ If $ \; r=s^{-1}$ where $s\in \p$,   then $r\in
  \p$,
  and $dr=|s|^{-m-1}ds$.
\end{lemma}

  For $Re
 \, \a >d-1$, the  Siegel gamma  function of $\p$
 is defined by
\be\label{2.444}
 \gm (\a)=\int_{\p} \exp(-\tr (r)) |r|^{\a } d_*r
 =\pi^{m(m-1)/4}\prod\limits_{j=0}^{m-1} \Gam (\a- j/2), \ee
\cite{FK, Te}. The relevant beta function has the form
\be\label{2.6}
 B_m (\a ,\b)=\int_0^{I_m} |r|^{\a -d} |I_m-r|^{\beta -d} dr=
 \frac{\gm (\a)\gm (\b)}{\gm (\a+\b)}, \quad d= (m+1)/2. \ee
This integral converges absolutely if and only if $Re
 \, \a, Re \, \b >d-1$.

 All function spaces on $\Ma$ are identified
with the corresponding spaces on $\bbr^{nm}$. For instance,
$\S(\Ma)$ denotes the Schwartz space    of infinitely differentiable
rapidly decreasing functions. The Fourier transform
 of a function $f\in L_{1}(\Ma)$ is defined by \be\label{ft}
\F f(y)=\int_{\Ma} \exp(\tr(iy'x)) f (x) dx,\qquad y\in\Ma \; .\ee

The {\it Cayley-Laplace operator} $\Del$
 on  $ \Ma$ is
defined by \be\label{K-L} \Del=\det(\d '\d), \qquad \d=(\d/\d
x_{i,j}). \ee In terms of the Fourier transform, the action of
$\Del$ represents a multiplication by the homogeneous polynomial
$(-1)^m |y|_m^2$  of degree $2m$ of $nm$ variables $y_{i,j}$.

For the sake of simplicity, for some operators on functions of
matrix argument  we will use the same notation as in the previous
sections.

The G{\aa}rding-Gindikin integrals of functions $f$ on $\p$
  are defined by
\be\label{3.1} (I_{+}^\a f)(s) \! =  \! \frac {1}{\gma} \intl_0^s \!
f(r)|s \! - \! r|^{\a-d} dr, \quad  (I_{-}^\a f)(s) \!  =  \! \frac
{1}{\gma} \intl_s^\infty  \! f(r)|r \! - \! s|^{\a-d} dr,\ee where
$s \in \p$ in the first integral and $s \in \cpm$ in the second one.
We assume $Re \, \a > d-1$, $ d=(m+1)/2$ (this condition is
necessary for absolute convergence of these integrals). The first
integral exists a.e. for arbitrary locally integrable function $f$.
Existence of the second integral requires extra assumptions for $f$
at infinity.

The {\it Riesz potential} of a function $f\in\S(\Ma)$
 is defined by \be\label{rie} (I^\a f)(x)=\frac{1}{\gam_{n,m} (\a)} \int_{\Ma}
f(x-y) |y|^{\a-n}_m dy;\ee \be\label{gam} \gam_{n,m}
(\a)=\frac{2^{\a m} \, \pi^{nm/2}\, \Gam_m (\a/2)}{\Gam_m
((n-\a)/2)}, \ \ \ Re \, \a>m-1, \ \ \a\neq n-m+1, \,  n-m+2, \ldots
\ee This integral is finite a.e. for $f \in L_{p}(\Ma)$ provided $ 1
\le p <n(Re \, \a +m-1)^{-1}$ \cite[Theorem 5.10]{R9}.

An application of the Fourier transform gives \be\label{hek} \F
[I^\a f](\xi)=|\xi|_m^{-\a} \F f(\xi)\ee (as in the case of $\rn$),
so that $I^\a$ can be formally identified with the negative power of
the  Cayley-Laplace operator (\ref{K-L}), namely,
$I^\a=(-\Del_m)^{-\a/2}$. Discussion of precise meaning of the
equality (\ref{hek}) and related references can be found in
\cite{R9}, \cite {OR2}.

\begin{definition} For $x \in \Ma, \;  n \ge m$, and $t \in \p$, we
define the (generalized) heat kernel $h_t (x)$ by the formula
\be\label{heat} h_t (x)=(4\pi)^{-nm/2}|t|^{-n/2} \exp (-\tr (t^{-1}
x'x)/4), \qquad |t|=\det (t), \ee and set \be\label{ga} H_t
f(x)=\intl_{\Ma} h_t (x-y) f(y) dy=\intl_{\Ma} h_{I_m} (y)
f(x-yt^{1/2}) \, dy. \ee
\end{definition}

 Clearly, $H_t f(x)$ is a  generalization of the
Gauss-Weierstrass integral (\ref{2.4}).

\begin{lemma}\label{hke}\cite{R9}{}\hfil

\noindent {\rm (i)} For each $ \; t \in \p$, \be\label{ed}
\int_{\Ma} h_t (x) \,dx =1.\ee

\noindent {\rm (ii)} The Fourier transform of $ \; h_t (x)$ has
the form \be\label{ft}\F h_t(y)= \exp (-\tr (ty'y), \ee which
implies the semi-group property  \be\label{cnv} h_t \ast
h_\tau=h_{t+\tau}, \qquad t, \tau \in \p. \ee

\noindent {\rm (iii)} If $f \in L_{p}(\Ma), \; 1\le p \le \infty$,
then \be\label{gw} ||H_t f||_p \le ||f||_p \, , \qquad \quad H_t
H_\tau f=H_{t+\tau}f, \ee and \be\label{lim}\lim\limits_{t \to
0}(H_t f)(x)=f(x) \ee
 in the $L_{p}$-norm. If $f$ is  a continuous function vanishing at infinity, then
 (\ref{lim}) holds  in the $\sup$-norm.
\end{lemma}
\begin{theorem}\label{rrg} \cite{R9}  \ Let $ \; m-1<Re \, \a<n-m+1$, $d=(m+1)/2$. Then
\be\label{rg} (I^\a f)(x) = \frac {1}{\Gam_m(\a/2)} \int_{\p}
|t|^{\a/2}H_t f(x) \,d_*t, \qquad  d_*t=|t|^{-d}\, dt,\ee
\be\label{rgg} H_t [I^\a f](x) =I_{-}^{\a/2}[H_{(\cdot)}
f(x)](t),\ee provided that integrals on either side of the
corresponding equality exist in the Lebesgue sense.
\end{theorem}

\subsection{Composite wavelet transforms: open problems}
 Formula (\ref{rg}) provokes a natural construction of the relevant
 composite wavelet transform on $\Ma$ associated with the heat
 kernel and containing a $\p$-valued scaling parameter.  To find
 this construction, we first obtain an auxiliary integral
 representation of a power function of the form $|t|^{\lam -d}$, $d=(m+1)/2$.
 \begin{definition} A function $w$ on $\p$ is said to be symmetric
 if \be\label{defs} w(g\eta g^{-1})= w (\eta) \quad \text{for all} \quad g \in
 GL(m,\bbr), \quad \eta\in \p.\ee
\end{definition}
Note that if $w$ is  symmetric, then for any $s,t \in \p$,
\be\label{wsy} w(t^{1/2} s t^{1/2})= w(s^{1/2} t s^{1/2})\quad
\text{and} \quad w(ts)=w(st).\ee Indeed, the second equality follows
from (\ref{defs}) if we set $\eta= ts, \; g=t^{-1}$. The first
equality in (\ref{wsy}) is a consequence of the second one:
 $$ w(t^{1/2} s t^{1/2})= w(t^{-1/2}[t^{1/2} s
 t^{1/2}]t^{1/2})=w(st)=w(ts)=w(s^{1/2} t s^{1/2}).$$
\begin{lemma} Let $w$ be a symmetric function on $\p$ satisfying \be
\int_{\p}\frac{|w(\eta)|}{|\eta|^\lam}\, d\eta < \infty, \qquad
c=\int_{\p}\frac{w(\eta)}{|\eta|^\lam}\, d\eta\neq 0, \qquad
|\eta|=\det (\eta). \ee Then for $t \in \p$, \be\label{gra}
|t|^{\lam -d}=c^{-1}\, \int_{\p}\frac{w(a^{-1}t)}{|a|^{m+1-\lam}}\,
da, \qquad d=(m+1)/2.\ee
\end{lemma}
\begin{proof} By (\ref{wsy}) we have (set $a=\rho^{-1}, \;da=
\rho^{-2d}d\rho$ ) \bea \int_{\p}\frac{w(a^{-1}t)}{|a|^{m+1-\lam}}\,
da&=&\int_{\p}\frac{w(t^{1/2}a^{-1}t^{1/2})}{|a|^{m+1-\lam}}\,
da=\int_{\p}\frac{w(t^{1/2}\rho t^{1/2})}{|\rho|^{\lam-d}}\, d_*\rho
\nonumber\\
&=&|t|^{\lam -d}\int_{\p}\frac{w(\eta)}{|\eta|^\lam}\, d\eta=c\,
|t|^{\lam -d}.\nonumber\eea
\end{proof}

Now we  replace a power function in (\ref{rg}) according to
(\ref{gra}) with $\lam=\a/2$. For $Re \, \a> (m-1)/2$, we obtain
\bea (I^\a f)(x) &=& \frac {c^{-1}}{\Gam_m(\a/2)} \int_{\p} H_t
f(x)\,
dt \int_{\p}\frac{w(a^{-1}t)}{|a|^{m+1-\a/2}}\, da\nonumber \\
&=&\frac {c^{-1}}{\Gam_m(\a/2)} \int_{\p}\frac{d_*
a}{|a|^{d-\a/2}}\int_{\p} H_t f(x) \,w (a^{-1}t)\, dt.\nonumber \eea
This gives \be \label {prot}(I^\a f)(x) = \frac
{c^{-1}}{\Gam_m(\a/2)}\int_{\p} \H f(x,a)|a|^{\a/2}  \,d_*a,\qquad
Re \, \a> (m-1)/2,\ee with \be \H f(x,a)=|a|^{-d}\int_{\p}  H_t
f(x)\,w (a^{-1}t)\, dt\ee or, by the symmetry of $w$, after changing
variable, \be \label {hewtr}\H f(x,a)=\int_{\p} H_{a^{1/2}\eta
a^{1/2}} f(x) w(\eta)\, d\eta, \qquad x \in \Ma, \quad a \in \p.\ee
Taking into account an obvious similarity between (\ref{hewtr})
 and the corresponding
``rank-one" formula for $m=1$, we call  $\H f(x,a)$ the {\it
composite wavelet transform of $f$ associated to the heat semigroup}
$H_t$. Here $w$ is a symmetric integrable function on $\p$ (that
will be endowed later with some cancelation properties) and $a$ is a
$\p$-valued scaling parameter.  One can replace $w$ by a more
general  {\it wavelet measure}, as we did in the previous sections,
but here we  want to  minimize technicalities.

Owing to (\ref{hek}), it is natural to expect that the inverse of
$I^\a$ has the same form (\ref{prot}) with $\a$ formally replaced by
$-\a$, and the case $\a=0$  gives a variant of Calder\'{o}n's
reproducing formula.

Thus, we encounter the following open problem:

{\bf Problem A.} {\it Give precise meaning to the inversion formula
\be\label{rghy} f(x) =c_{m,\a} \int_{\p} \frac{\H \vp
(x,a)}{|a|^{\a/2}} \,d_*a, \qquad \vp=I^\a f,\ee and the reproducing
formula \be\label{rghr} f(x) =c_{m} \int_{\p} \H f (x,a)\, d_*a,\ee
say, for $ f\in L_p$ or any other ``natural" function space. Give
examples of wavelet functions $w$ for which (\ref{rghy}) and
(\ref{rghr}) hold. Find explicit formulas for the normalizing
coefficients $c_{m,\a}$ and $c_{m}$, depending on $w$.}

Solution of this problem would give a series of pointwise inversion
formulas for diverse Radon-like transforms on matrix spaces; see,
e.g., \cite {OR2}, \cite {OR3}, \cite {R9},
 where such formulas are available in terms of distributions.
 Justification of (\ref{rghy}) and (\ref{rghr}) would also bring new light
 to a variety of inversion formulas for Radon transforms on
 Grassmannians, cf. \cite {GRu}.

\subsection{Some discussion} Trying to solve Problem A, we  come across new
problems that are of independent interest.
 Let $Re \, \a>d-1, \; d=(m+1)/2$. Suppose, for instance, that $f(x), \; x \in \Ma$, is
 a Schwartz function and $w (\eta), \; \eta \in \p$, is ``good
 enough". We anticipate the following equality:
 \be\label{ant} I_\e f(x)\equiv \int_{\e I_m}^\infty\frac{\H [I^\a
 f]
(x,a)}{|a|^{\a/2}} \,d_*a= \int_{\p} \Lam_{\a/2}(s)\, H_{\e s}
f(x)\, ds,\ee where $\Lam_{\a/2}(s)$ expresses through the
G{\aa}rding-Gindikin integral in (\ref{3.1}) as \be\label {ky}
\Lam_{\a/2}(s)=\frac{\Gam_m (d)}{|s|^d}\, I_+^{\a/2 +d} w (s),
\qquad s\in \p.\ee If $m=1$ and $\a/2$ is replaced by $\a$, then
(\ref{ky}) coincides with the function $\lambda _{\alpha }(s)=s^{-1}
I_+^{\alpha +1}\mu (s)$ in Lemma \ref{lB}. Now, we give the
following
\begin{definition} An integrable symmetric  function $w$ on $\p$ is
called an {\it admissible wavelet} if \be\label {wad}
\Lam_{\a/2}(s)\equiv \frac{\Gam_m (d)}{|s|^d}\, I_+^{\a/2 +d} w (s)
\in L_1 (\p)\quad \text{\rm and} \quad
c_\a=\int_{\p}\Lam_{\a/2}(s)\,ds\neq 0.\ee
\end{definition}
If $w$ is admissible, then, by Lemma \ref {hke}, the $L_p$-limit as
$\e \to 0$ of the right-hand side of (\ref{ant}) is $c_\a \,f$, and
we are done. This discussion includes the case $\a=0$ corresponding
to the reproducing formula.

Thus, our attempt to solve Problem A rests upon the following

{\bf Problem B.} {\it Find examples of admissible wavelets (both for
$\a\neq 0$ and $\a=0$) and compute $c_\a$.}

Now, let us try to prove (\ref{ant}). We say ``try", because along
the way, we  come across one more open problem related to
application of the Fubini theorem; cf. justification of interchange
of the order of integration in the proof of Theorem \ref{teo:34}.

By (\ref{hewtr}) and (\ref{rgg}), \bea \H I^\a f(x,a)&=&\int_{\p}
H_{a^{1/2}\eta a^{1/2}} I^\a f(x) \,w(\eta)\, d\eta\nonumber \\
&=&\int_{\p} I_{-}^{\a/2}[H_{(\cdot)} f(x)](a^{1/2}\eta
a^{1/2})\,w(\eta)\, d\eta.\nonumber \eea Assume that $x$ is fixed
and denote $\psi (s)=H_{s} f(x)$. Then \bea
 \H I^\a f(x,a)&=&\int_{\p}w(\eta)\,I_{-}^{\a/2}\psi (a^{1/2}\eta a^{1/2})\,
 d\eta\nonumber \\
&=& \frac {1}{\Gam_m(\a/2)} \int_{\p}w(\eta)\,d\eta
\int_{a^{1/2}\eta
 a^{1/2}}^\infty \psi (s) |s-a^{1/2}\eta a^{1/2}|^{\a/2 -d}\,
 ds\nonumber \\
&=&\frac {1}{\Gam_m(\a/2)} \int_{\p}\psi (s)\,
ds\int_0^{a^{-1/2}\eta a^{-1/2}}w(\eta)\,|s-a^{1/2}\eta
a^{1/2}|^{\a/2 -d}\,d\eta\nonumber \\
&=&|a|^{\a/2 -d} \int_{\p}\psi (s)\,I_{+}^{\a/2}w (a^{-1/2}s
a^{-1/2})\, ds. \nonumber \eea Hence, the left-hand side of
(\ref{ant}) transforms as follows.

\bea I_\e f(x)&=& \int_{\e I_m}^\infty \frac{da}{|a|^{m+1}}
\int_{\p}\psi (s)\,I_{+}^{\a/2}w (a^{-1/2}s a^{-1/2})\, ds\nonumber \\
&=&\int_{\p}\psi (s)\,\int_{\e I_m}^\infty I_{+}^{\a/2}w (a^{-1/2}s
a^{-1/2})\frac{da}{|a|^{m+1}} \qquad \text{\rm (set
$a=\t^{-1}$)}\nonumber \\
&=&\int_{\p}\psi (s)\,\int_0^{\e^{-1}I_m} I_{+}^{\a/2}w (\t^{1/2}s
\t^{1/2})\, d\t\nonumber \\
&=&\e^{md}\int_{\p}\psi (\e s)\, dy\int_0^{\e^{-1}I_m} I_{+}^{\a/2}w
(\t^{1/2}\e^{1/2}s\e^{1/2} \t^{1/2})\, d\t.\nonumber\eea Thus we
have \be I_\e f(x)=\int_{\p}\psi (\e s)\,k(s)\, ds=\int_{\p}H_{\e s}
f(x) \,k(s)\, ds, \ee where $$k(s)=\int_0^{I_m}
I_{+}^{\a/2}w(\lam^{1/2}s\lam^{1/2})\, d\lam.$$ To get (\ref{ant}),
it remains to show that $k(s)$ coincides with the function
(\ref{ky}). We have \bea k(s)&=&\frac
{1}{\Gam_m(\a/2)}\int_0^{I_m}d\lam \int_0^{\lam^{1/2}s\lam^{1/2}}
w(s) |\lam^{1/2}s\lam^{1/2} -s|^{\a/2 -d}\, ds\nonumber \\
&&\text{(set $s=\lam^{1/2}z\lam^{1/2}$ and note that
$w(\lam^{1/2}z\lam^{1/2})=w(z^{1/2}\lam z^{1/2})$)}\nonumber \\
&=&\frac {1}{\Gam_m(\a/2)}\int_0^{I_m} |\lam|^{\a/2}d\lam \int_0^s
|s-z|^{\a/2 -d}w(z^{1/2}\lam z^{1/2})\, dz\nonumber \\
&=&\frac {1}{\Gam_m(\a/2)} \int_0^s |s-z|^{\a/2 -d}\, dz\int_0^{I_m}
|\lam|^{\a/2}\,w(z^{1/2}\lam z^{1/2})\,d\lam\nonumber \\
&=&\frac {1}{\Gam_m(\a/2)} \int_0^s |s-z|^{\a/2 -d}\,
\frac{dz}{|z|^{\a/2 +d}}\int_0^z w(b) |b|^{\a/2}\, db\nonumber \\
&=&\frac {1}{\Gam_m(\a/2)} \int_0^s  w(b) |b|^{\a/2}\,u(b,s)\,
db,\nonumber \eea where
 \bea u(b,s)&=& \int_b^s |s-z|^{\a/2 -d}\,
\frac{dz}{|z|^{\a/2 +d}}
\qquad \text{(set $z=r^{-1}$)}\nonumber \\
&=& \int_{s^{-1}}^{b^{-1}}|sr -I_m|^{\a/2 -d}\, dr=|s|^{\a/2
-d}\int_{s^{-1}}^{b^{-1}}|r -s^{-1}|^{\a/2 -d}\, dr.\nonumber \eea
The last integral can be easily computed using the well-known
formula for Siegel Beta functions \be\intl_a^b |r-a|^{\a -d}
|b-r|^{\beta -d} dr= B_m (\a ,\b) |b-a|^{\a+\beta -d}\ee (many such
formulas can be found, e.g., in \cite{OR2}), and we have \be u(b,s)=
B_m (\a/2, d)\, \frac{|s-b|^{\a/2}}{|s|^{d}\, |b|^{\a/2}}, \qquad
B_m (\a/2, d)=\frac{\Gam_m(\a/2)\, \Gam_m(d)}{\Gam_m(\a/2 +d)}.\ee
Finally, we  get
$$ k(s)=\frac{\Gam_m(d)}{|s|^{d}\,\Gam_m(\a/2 +d)} \int_0^s
w(b)|s-b|^{\a/2}\, ds=\frac{\Gam_m(d)}{|s|^{d}}\, I^{\a/2 +d}_+
w(s)= \Lam_{\a/2}(s).$$

{\bf Problem C.} Although all calculations above go through
smoothly, interchange of the order of integration remains
unjustified. We do not know how to justify it and what  additional
requirements  on the  wavelet $w$  should be imposed (if any). One
of the obstacles  is that $\int_0^\infty \neq \int_0^s
+\int_s^\infty$, when we integrate over the higher-rank  cone.

\end{document}